\documentclass[a4paper,12pt]{article}
\usepackage{a4wide}
\usepackage{amsthm,amsmath}
\usepackage{amssymb}
\usepackage[all]{xypic}
\usepackage{enumerate}

\begin{document}

\swapnumbers
\newtheorem{theorem}{Theorem}[section]
\newtheorem{proposition}[theorem]{Proposition}
\newtheorem{corollary}[theorem]{Corollary}
\newtheorem{lemma}[theorem]{Lemma}
\newtheorem*{theorem*}{Theorem}
\newtheorem*{proposition*}{Proposition}
\newtheorem*{corollary*}{Corollary}
\newtheorem*{lemma*}{Lemma}\theoremstyle{definition}
\newtheorem{definition}[theorem]{Definition}
\newtheorem{punto}[theorem]{}
\newtheorem{example}[theorem]{Example}
\newtheorem{remark}[theorem]{Remark}
\newtheorem{remarks}[theorem]{Remarks}
\newtheorem*{definition*}{Definition}
\newtheorem*{remark*}{Remark}

\newcommand{\cat}[1]{\mathcal{#1}}
\newcommand{\unidad}[1]{\eta_{#1}}
\newcommand{\counidad}[1]{\delta_{#1}}
\newcommand{\stat}[1]{\mathrm{Stat}(#1)}
\newcommand{\adstat}[1]{\mathrm{Adst}(#1)}
\newcommand{\coker}[1]{\mathrm{coker}(#1)}
\newcommand{\rcomod}[1]{\mathcal{M}^{#1}}
\newcommand{\lcomod}[1]{{}^{#1}\mathcal{M}}
\newcommand{\rmod}[1]{\mathcal{M}_{#1}}
\newcommand{\lmod}[1]{{}_{#1}\mathcal{M}}
\newcommand{\coend}[2]{\mathrm{Coend}_{#1}(#2)}
\newcommand{\End}[3]{\mathrm{End}_{#1}^{#2}(#3)}
\renewcommand{\hom}[3]{\mathrm{Hom}_{#1}(#2,#3)}
\newcommand{\Hom}{\mathrm{Hom}}
\newcommand{\cotensor}[1]{\square_{#1}}
\newcommand{\functor}[1]{\mathbf{#1}}
\newcommand{\Mod}{\mbox{{\rm -Mod}}}
\newcommand{\cohom}[3]{\mathrm{h}_{#1}(#2,#3)}
\newcommand{\HOM}[3]{\mathrm{HOM}_{#1}(#2,#3)}

%

\title{On quasi-Frobenius bimodules and corings}
\author{F.~Casta\~no Iglesias
  \thanks{Research partially
supported by Spanish Project (MTM2005--03227) from MCT.}
 \\{\small Departamento de
Estad\'{\i}stica y Matem\'atica Aplicada,} \\{\small Universidad de
Almer\'{\i}a 04120, Almer\'{\i}a, Spain}\\{\small E-mail:
fci@ual.es}}
\date{}
\maketitle
%
%

\section{Introduction}
Frobenius bimodules  are connected with Frobenius algebras and
extensions. For instance, a ring extension $\varphi:R\rightarrow S$
is a Frobenius extension if and only if $_RS_S$ is a Frobenius
bimodule \cite{A-F}. Brzezi\'{n}ski and G\'omez-Torrecillas studied
in \cite{BG} certain properties of comatrix corings in relation to
properties
 of bimodules. In particular  they showed that the \emph{comatrix coring}
\cite{kaouti-pepe} induced by any Frobenius bimodule is a Frobenius
coring. Recently, Casta\~no Iglesias and N{\u a}st{\u a}sescu
\cite{CN} studied quasi-Frobenius extensions  from the viewpoint of
quasi-Frobenius functors. The ring extension $\varphi:R\rightarrow
S$ is quasi-Frobenius (left and right extension) if and only if the
restriction of scalars functor $\varphi_*:\lmod{S}\rightarrow
\lmod{R}$ is a quasi-Frobenius functor (Remark \ref{restri}).
\\
Following this direction, we shall extend some results in \cite{BG}
to the case when $_RM_S$ is a \emph{quasi-Frobenius bimodule}. Using
the notion of quasi-Frobenius functor and properties on
quasi-Frobenius corings given in \cite{CN}, we discuss the questions
of when a quasi-Frobenius bimodule induce a quasi-Frobenius
extension and a quasi-Frobenius coring.
\bigskip

Throughout this paper all rings will be assumed to have identity,
and all modules to be unital. We use the standard module theory
notation, for example a $(R,S)$-bimodule $M$ is denoted by
$_RM_S$, $\hom{R}{-}{-}$ denotes the Abelian group of $R$-module
maps. The dual of a left $R$-module $M$ is denoted by $(_RM)^*$.
Finally, $_R\mathcal{M}_S$ denotes the category of all
$(R,S)$-bimodules. Let
 $_RM_S$ and $_RN_S$ be $(R,S)$-bimodules. $M$ is said to be \emph{similar} to $N$, abbreviated $_RM_S\sim \, _RN_S$,
 if there
are $m,n\in \mathbb{N}$ and $(R,S)$-bimodules $P$ and $Q$ such
that $M\oplus P \cong N^{(m)}$ and $N\oplus Q \cong M^{(n)}$ as
$(R,S)$-bimodules (cf. \cite{A-F}). It is easy see that ``$\sim
$'' defines an equivalence relation on the class  of
$(R,S)$-bimodules. We start by recalling the definition of
quasi-Frobenius functors between Grothendieck categories given in
\cite{CN}.

 Let $\mathcal{A}$ be a
Grothendieck category and consider an object $X$ of $\cat{A}$. For
any positive integer $n$, we denote by $X^n$ the direct sum of $n$
copies of $X$ in the category $\cat{A}.$ Given additive and
covariant functors $\functor{L,R}: \mathcal{B}\rightarrow
\mathcal{A}$ between Grothendieck categories, we say that
$\functor{L}$ \emph{divide} to $\functor{R}$, denoted by
$\functor{L}\mid\functor{R},$  if for some positive integer $n$
there are natural transformations
\[\xymatrix{\functor{L}(X)\ar[r]^-{\phi(X)}&  \functor{R}(X)^n
 \ar[r]^-{\psi(X)}& \functor{L}(X) }\]
such that $\psi(X)\circ\phi(X) = 1_{\functor{L}(X)}$ for every $X\in
\mathcal{B}$.  $\functor{R}\mid\functor{L}$ is defined
symmetrically. The functor $\functor{L}$ is said to be
 \emph{similar} to $\functor{R}$, denoted by
 $\functor{L}\sim\functor{R},$ if  $\functor{L}\mid\functor{R}$ and
$\functor{R}\mid\functor{L}.$

\begin{definition}
The functor $ \functor{F}:\mathcal{A}\rightarrow \mathcal{B}$ is
said to be \emph{quasi-Frobenius functor} if $\functor{F}$ has a
left adjoint $\functor{L}: \mathcal{B}\rightarrow \mathcal{A}$ and
also a right adjoint $\functor{R}:\mathcal{B}\rightarrow
\mathcal{A}$ with $\functor{L}$ similar to $\functor{R}$.
\end{definition}

\begin{remark}\label{restri}
 Let $\varphi: R\rightarrow S$ be a ring extension. Associated  to $\varphi$ we have the
restriction of scalars functor $\varphi_*: \lmod{S}\rightarrow
\lmod{R}$, the induction functor $S\otimes_R-: \lmod{R}\rightarrow
\lmod{S}$ and the coinduction functor
$\hom{R}{_RS_S}{-}:\lmod{R}\rightarrow \lmod{S}$. It is well known
that $\varphi_*$ is right adjoint to $S\otimes_R-$  and left adjoint
to $\hom{R}{_RS_S}{-}$. Then $\varphi_*$ is a \emph{quasi-Frobenius
functor} if the functor $S\otimes_R-$ is similar to
$\hom{R}{_RS_S}{-}$. This is equivalent to assuming $_RS$ finitely
generated projective and $_SS_R\sim \, \hom{R}{_RS}{R}$. On the
other hand, recall of \cite{muller} that a ring homomorphism
$\varphi: R\rightarrow S$ is called \emph{left quasi-Frobenius
extension} if $_RS$ is finitely generated and projective and  $_SS_R
\, | \, \hom{R}{_RS}{R}$. Equivalently, $S_R$ and $_RS$ are finitely
generated  projective and  $ \hom{R}{S_R}{R} \, | \, _RS_S$.
\\
Similarly, $\varphi$ is called  \emph{right quasi-Frobenius
extension} if $S_R$ is finitely generated projective and
 $_RS_S \, | \, \hom{R}{S_R}{R}$ (or  $
\hom{R}{_RS}{R} \, | \,  _SS_R$).
 Then $\varphi$ is
a \emph{quasi-Frobenius extension} (left and right extension) if
$_RS$  is a finitely generated projective  and $_SS_R \sim
\hom{R}{_RS}{R}$ (or $_RS_S \sim \hom{R}{S_R}{R}$). Therefore, the
ring extension $\varphi$ is quasi-Frobenius if and only if
$\varphi_*:\lmod{S}\rightarrow \lmod{R}$ is a quasi-Frobenius
functor.
\end{remark}

\section{Quasi-Frobenius bimodules}
In this section we generalize  the notion of Frobenius bimodule
 to \emph{quasi-Frobenius bimodule}. Recall of \cite{A-F} that a $(R,S)$-bimodule $M$ is called Frobenius bimodule if $_RM$ and $M_S$ are
finitely generated and projective modules and $(_RM)^*\cong \,
(M_S)^*$ as $(S,R)$-bimodules. We start with the following
definition.
\begin{definition}\label{qF-bimodule} An  $(R,S)$-bimodule $M$ is
said to be \emph{quasi-Frobenius bimodule}, if both $_RM$ and $M_S$
are finitely generated projective and $(_RM)^*\sim \, (M_S)^*$ as
$(S,R)$-bimodules.
\end{definition}
For any ring extension $\varphi:R\rightarrow S,$ it is easy to see
that $\varphi$ is a quasi-Frobenius extension if and only if the
natural bimodule $_RS_S$ is a quasi-Frobenius bimodule. So
quasi-Frobenius bimodules generalize quasi-Frobenius extension.
\begin{remark} It is obvious that the class of quasi-Frobenius
bimodule contains Frobenius bimodules. In general, a quasi-Frobenius
bimodule need not be a Frobenius (an example of finite-dimensional
quasi-Frobenius algebra which is not a Frobenius is gives in
\cite{naka}).
\end{remark}

The following shows some relations between quasi-Frobenius
functors and bimodules. If $_RM_S$ and $_SN_T$ are bimodules, then
$M\otimes_SN$ receives the natural $(R,T)$-bimodule structure
indicated by $r(m\otimes n)t=rm\otimes nt$.
\begin{proposition}\label{product} Suppose $_RM_S$ and $_SN_T$ are quasi-Frobenius bimodules. Then $_R(M\otimes_SN)_T$ is  a
quasi-Frobenius $(R,T)$-bimodule.
\end{proposition}
\begin{proof} First, note that $(_R(M\otimes_SN))^*$ and $((M\otimes_SN)_T)^*$ are  finitely
generated and projective modules, since $_RM, \, M_S, \, _SN$ and
$N_T$ are finitely generated projective modules. Next we apply
$\hom{S}{_SN}{-}$ to $(_RM)^*\sim \, (M_S)^*$ to obtain:
$$\hom{S}{_SN}{(_RM)^*}\sim \, \hom{S}{_SN}{(M_S)^*}$$
Equivalently, $(_R(M\otimes_SN))^* \sim \, ((M\otimes_SN)_S)^*$.
Similarly, applying $\hom{S}{M_S}{-}$ to $(_SN)^*\sim \, (N_T)^*$,
it follows that $((M\otimes_SN)_S)^*\sim \,((M\otimes_SN)_T)^*$.
 Thus $(_R(M\otimes_SN))^* \sim ((M\otimes_SN)_T)^*,$ whence $M\otimes_SN$
is a quasi-Frobenius $(R,T)$-bimodule.
\end{proof}

\begin{corollary} If $\varphi:R\rightarrow S$  and $\psi:S\rightarrow T$
are  quasi-Frobenius ring extensions, then $_RT_T$ is a
quasi-Frobenius bimodule.
 \end{corollary}
\begin{proof} The corollary follows from letting $M=\, _RS_S$ and $N=\, _ST_T$
in the proposition.
\end{proof}
Let $M$ be any $(R,S)$-bimodule. It is well-known that the functor
$M\otimes_S-:\lmod{S}\rightarrow \lmod{R}$ has a right adjoint
$\hom{R}{_RM}{-}: \lmod{R}\rightarrow \lmod{S}.$ If in addition,
$M_S$ is a finitely generated and projective module, then the
functor $M\otimes_S-$ has also a left adjoint $(M_S)^*\otimes_R- :
\lmod{R}\rightarrow \lmod{S}.$ When this is the case,
$M\otimes_S-$ is a quasi-Frobenius functor if the functors
$(M_S)^*\otimes_R-$ and $\hom{R}{_RM}{-}$ are similar. This leads
to a proposition needed later:
\begin{proposition}\label{bimo-funtor} For any $(R,S)$-bimodule $M$ the following
assertions are equivalent.
 \begin{enumerate}[(a)]
 \item $_RM_S$ is a quasi-Frobenius bimodule;
 \item $M\otimes_S-$
is a quasi-Frobenius functor.
\end{enumerate}
\end{proposition}
\begin{proof} $(a)\Rightarrow (b).$ Suppose $_RM$ and $M_S$ are finitely generated projective modules and $(M_S)^* \,
\sim \, (_RM)^*.$ Then  by \cite[Lemma 2.1]{CN}, $(M_S)^*\otimes_R-
\, \sim \, (_RM)^*\otimes_R-$. We note that $(_RM)^*\otimes_R-\cong
\hom{R}{_RM}{-}$, since $_RM$ is a finitely generated projective.
Therefore, $(M_S)^*\otimes_R-$ and $\hom{R}{_RM}{-}$ are similar,
which implies that  $M\otimes_S-$ is a quasi-Frobenius functor.
\\
$(b)\Rightarrow (a).$ Assume that $M\otimes_S-$ is a quasi-Frobenius
functor. By \cite[Theorem 2.2]{CN}, $_RM_S$ is a finitely generated
 projective  on both sides with $(M_S)^*\otimes_R- \sim \,
\hom{R}{_RM}{-}.$ So $(M_S)^*\otimes_RR \sim \, \hom{R}{_RM}{R},$
whence $(M_S)^* \, \sim \, (_RM)^*$. Consequently, $_RM_S$ is a
quasi-Frobenius bimodule.
\end{proof}
As a immediate consequence of the proposition above with $M= \,
_RS_S$ and Remark \ref{restri}, we have
\begin{corollary} Let $\varphi:R\rightarrow S$ be a ring extension.
Are equivalent the following assertions
\begin{enumerate}[(i)]
\item
 $\varphi$ is a
quasi-Frobenius extension. \item $_RS_S$ is a quasi-Frobenius
bimodule. \item  $\varphi_*: \lmod{S}\rightarrow \lmod{R}$ is a
quasi-Frobenius functor.
\end{enumerate}
\end{corollary}

Consider now an $(R,S)$-bimodule $M$ and let $T=\End{S}{}{M_S}$ be
the endomorphism ring of $M_S$. Then there is a ring extension
$\overline{\varphi}:R\rightarrow T$ defined by $r\mapsto
\overline{\varphi}_r(m)=rm$. Note too the natural bimodule $_TM_S$
given by $tms=t(m)s$. The next theorem generalize the classical
endomorphism ring theorem to the case of quasi-Frobenius bimodules.
\begin{theorem}
If $_RM_S$ is a quasi-Frobenius bimodule, then $T=\End{S}{}{M_S}$ is
a quasi-Frobenius extension of $R$.
\end{theorem}
\begin{proof} By Definition \ref{qF-bimodule}, it suffices to prove that $_RT_T$ is a quasi-Frobenius
bimodule. Applying  $M\otimes_S-$ to $(M_S)^* \, \sim \, (_RM)^*,$
we obtain $M\otimes_S(M_S)^* \, \sim \, M\otimes_S(_RM)^*.$ This
establishes that $ T= \End{S}{}{M_S}\cong M\otimes_S(M_S)^* \, \sim
\, M\otimes_S(_RM)^* = (_RT)^**$. On the other hand, note that $_RT$
is finitely generated projective, since $M_S$ and $(M_S)$ are
 finitely generated projective modules. So $_RT_T$ is a quasi-Frobenius
bimodule.
\end{proof}
A converse of this theorem is give below where (\ref{willard}) is a
\emph{condition Willard's}, which is satisfied by a generator module
\cite{willa}.
\begin{proposition} Let $R$ and $S$ be rings and $M$  a
$(R,S)$-bimodule such that $_RM$ and $M_S$ are both finitely
generated projective modules. Suppose $T=\End{S}{}{M_S}$ is a
quasi-Frobenius extension of $R$, and
\begin{equation}\label{willard}
\hom{T}{_TM_S}{_TT_R}\cong \hom{S}{_RM_S}{S_S}.
\end{equation}
Then $_RM_S$ is a quasi-Frobenius bimodule.
\end{proposition}
\begin{proof} In the computation below, we apply  the Hom-Tensor
Relation, the necessary condition $(_RT)^*\sim (T_T)^*$ for a
quasi-Frobenius extension $R\rightarrow T$ and condition Willard
$(_TM)^*\cong (M_S)^*$  in that order:
$$\begin{array}{ccl} (_RM)^* = \hom{R}{_RM}{R}& = &
          \hom{R}{T\otimes_TM}{R} \\[+1mm]
&\cong & \hom{T}{_TM}{\hom{R}{_RT}{R}} \\[+1mm]
&= & \hom{T}{_TM}{(_RT)^*}\\[+1mm]
 &\sim& \hom{T}{_TM}{(T_T)^*}\\[+1mm]
 &\cong& \hom{T}{_TM}{T} = (_TM)^* \cong (M_S)^*
\end{array}$$
Now the proposition follows, since $_RM$ and $M_S$ are both finitely
generated projective modules.
\end{proof}
\section{Application to corings}
 Let $A$ be an
associative and unitary algebra over a commutative ring (with unit)
$k$. We recall from \cite{Sw2} that an $A$-coring $\mathfrak{C}$
consists of an $A$-bimodule $\mathfrak{C}$ with two $A$-bimodule
maps
$$\Delta: \mathfrak{C}\rightarrow \mathfrak{C}\otimes_A\mathfrak{C},
\mbox{      } \epsilon: \mathfrak{C}\rightarrow A$$ such that
$(\mathfrak{C}\otimes_A\Delta)\circ \Delta =
(\Delta\otimes_A\mathfrak{C})\circ \Delta$ and
$(\mathfrak{C}\otimes_A\epsilon)\circ \Delta =
(\epsilon\otimes_A\mathfrak{C})\circ \Delta = 1_{\mathfrak{C}}.$ A
right $\mathfrak{C}$-comodule is a pair $(M,\rho_M)$ consisting of a
right $A$-module $M$ and an $A$-linear map $\rho_M:M\rightarrow
M\otimes_A\mathfrak{C}$ satisfying $(M\otimes_A\Delta)\circ \rho_M =
(\rho_M\otimes_A\mathfrak{C})\circ \rho_M$ and
$(M\otimes_A\epsilon)\circ \rho_M =1_M$. The right
$\mathfrak{C}$-comodules together  with their morphisms form the
additive category $\mathcal{M}^{\mathfrak{C}}$. If $_A\mathfrak{C}$
is flat, then $\mathcal{M}^{\mathfrak{C}}$ is a Grothendieck
category.

 The forgetful functor
$\functor{U}:\rcomod{\mathfrak{C}}\rightarrow \rmod{A}$ has the
right adjoint $-\otimes_A\mathfrak{C}$ (see \cite[Lemma 3.1]{Brze}).
When $_A\mathfrak{C}$ is finitely generated and projective, the
functor $\hom{A}{\mathfrak{C}}{-}:\rmod{A}\rightarrow
\rcomod{\mathfrak{C}}$ is a left adjoint to $\functor{U}$. In this
case, $\functor{U}$ is a quasi-Frobenius functor if
$-\otimes_A\mathfrak{C}$ and $\hom{A}{\mathfrak{C}}{-}$ are similar.
\begin{definition}  An A-coring $\mathfrak{C}$ is called
\emph{quasi-Frobenius coring} provided the forgetful functor
$\functor{U}:\rcomod{\mathfrak{C}}\rightarrow \rmod{A}$ is a
quasi-Frobenius functor.
\end{definition}
A characterization of such corings is given in \cite[Theorem
5.5]{CN}. In particular an $A$-coring $\mathfrak{C}$ is a
quasi-Frobenius coring if and only if $_A\mathfrak{C}$ is finitely
generated projective and $\mathfrak{C}\sim E$ as $(A,E)$-bimodules
where $E$ is the opposite algebra of $\mathfrak{C}^*$.
\\
Associated to any $(R,S)$-bimodule $M$ we can consider the
\emph{comatrix $S$-coring} $(M_S)^*\otimes_RM$ provided $M_S$ is a
finitely generated projective \cite{kaouti-pepe}. As for Frobenius
bimodules we obtain the following theorem, which generalizes
\cite[Theorem 3.7]{BG}.

\begin{theorem}\label{qF-coring} Let $_RM_S$ be a bimodule and $T=\End{S}{}{M_S}$ be the endomorphism ring. Then
\begin{enumerate}[(a)]
\item If $_RM_S$ is a quasi-Frobenius bimodule, then
$(M_S)^*\otimes_RM$ is a quasi-Frobenius $S$-coring. \item If
$(M_S)^*\otimes_RM$ is a quasi-Frobenius $S$-coring, then
$T\otimes_RT$ is a quasi-Frobenius $T$-coring.
\end{enumerate}
\end{theorem}
\begin{proof}$(a).$ Suppose that $_RM_S$ is a quasi-Frobenius bimodule. Then the functors $(M_S)^*\otimes_R-$ and
$\hom{R}{_RM}{-}$ are similar by Proposition \ref{bimo-funtor}. So
$(M_S)^*\otimes_RM \sim \, \hom{R}{_RM}{M}.$ From \cite[Proposition
2.1]{kaouti-pepe}, $\hom{R}{_RM}{M}\cong E$ where $E$ is the
opposite algebra of  $((M_S)^*\otimes_RM)^*$. Thus
$(M_S)^*\otimes_RM\sim E$. Now \cite[Theorem 4.4]{CN} implies that
$(M_S)^*\otimes_RM$ is a quasi-Frobenius coring, since
$(M_S)^*\otimes_RM$ is a finitely generated projective left
$S$-module by assumptions.
\\
$(b).$ Suppose that $(M_S)^*\otimes_RM$ is a quasi-Frobenius
$S$-coring. Then $(M_S)^*\otimes_RM \sim \, \hom{R}{_RM}{M}.$
Applying first the functor $M\otimes_S-$ and later
$-\otimes_S(M_S)^*$, we obtain
$$ M\otimes_S(M_S)^*\otimes_RM \otimes_S(M_S)^*\sim \, M\otimes_S\hom{R}{_RM}{M}\otimes_S(M_S)^* $$
This is equivalent to
$$T\otimes_RT \sim \, \hom{R}{T}{M\otimes_S(M_S)^*}\cong
\hom{R}{T}{T}\cong \overline{E}$$ where $\overline{E}$ is the
opposite algebra of $((_RT)^*\otimes_RT)^*$. Since
$(_RT)^*\otimes_RT \cong T\otimes_RT$, we note that  \cite[Theorem
4.4]{CN} finish the proof.
\end{proof}

 Following Sweedler \cite{Sw2}, given a ring extension
  $\rho:R\rightarrow S$ one can view $\mathfrak{C} = S\otimes_RS$ as
  an $S$-coring. $\mathfrak{C}$ is known as \emph{Sweedler's coring}
  associated to $\rho$. As an immediate consequence of Theorem \ref{qF-coring} we
  obtain the following corollary.
\begin{corollary}\label{sw-coring}
If $_RS_S$ is a quasi-Frobenius bimodule, then the Sweedler's
coring $S\otimes_RS$ is a quasi-Frobenius coring.
\end{corollary}
\begin{remark} A similar proof to Brzezi\'{n}ski \cite[Theorem
2.7]{Brze1} show that $R\rightarrow S$ is a quasi-Frobenius
extension whenever $_RS$ is faithfully flat and $S\otimes_RS$ is a
quasi-Frobenius coring.
\end{remark}

\begin{remark} Let $_RM_S$ be a bimodule and $T=\End{S}{}{M_S}$ be the ring endomorphism. Then the following diagram summarise the main results of
this paper.
\[
\xymatrix@*+<14pt>{*\txt{$_RM_S$ qF-bimodule  }
 \ar@{=>}[rrr]^--{\small{\txt{Th.
3.2(a)}}}\ar@{=>}[d]^{\small{\txt{Prop. 2.7}}} &&&
  *\txt{$(M_S)^*\otimes_RM$ qF-coring} \ar@{=>}[ddlll]^{\small{\txt{Th. 3.2(b)}}}  \\
 \txt{$R\rightarrow T$ qF-extension}\ar@<1ex>[u]^{\small{\txt{Willard's condition}}}_{}\ar@{=>}[d]^{\small{\txt{Cor. 3.3}}}&  \\
  \txt{$T\otimes_RT$ qF-coring} \ar@<1ex>[u]^{\small{\txt{Rem. 3.4}}}_{}& }
\]
\end{remark}

\bibliographystyle{amsplain}
\bibliography{e:/texdoc/bib/biblo1}

\end{document}